\documentclass[12pt]{amsart} 

\usepackage[utf8]{inputenc} 

\usepackage{lipsum}

\usepackage{graphicx} 

\usepackage[english]{babel}

\usepackage{booktabs} 
\usepackage{array} 
\usepackage{paralist} 
\usepackage{verbatim} 
\usepackage{subfig} 

\usepackage{caption}

\usepackage{amsfonts,amsmath,amssymb,bbm,amsthm,amsbsy}
\usepackage[alphabetic]{amsrefs}

\usepackage{mathtools}

\usepackage{color}
\usepackage{pdfcolmk}

\newcounter{dawidcomments}
\newcommand{\dawid}[1]{\textbf{\color{blue}(D\arabic{dawidcomments})} \marginpar{\scriptsize\raggedright\textbf{\color{blue}(D\arabic{dawidcomments})Dawid: }#1}
\addtocounter{dawidcomments}{1}}

\newcounter{ncomments}
\newcommand{\ncomm}[1]{\textbf{\color{red}(N\arabic{ncomments})} \marginpar{\scriptsize\raggedright\textbf{\color{red}(N\arabic{ncomments})Nici: }#1}
\addtocounter{ncomments}{1}}

\usepackage{tikz}
\usepackage{tikz-cd}

\usepackage{IEEEtrantools}

\newenvironment{equ*}[1]{\begin{IEEEeqnarray*}{#1}}{\end{IEEEeqnarray*}}

\usepackage[pdfborder={0 0 0 [0 0 ]},bookmarksdepth=3,bookmarksopen=true]{hyperref}
\usepackage[capitalize]{cleveref}
\usepackage[ps,all,arc,rotate]{xy}
\usepackage{bbm} 

\hypersetup{
    colorlinks=true,
    linkcolor=black,
    citecolor=black,
    filecolor=black,
    urlcolor=black,
}

\newenvironment{frcseries}{\fontfamily{frc}\selectfont}{}

\makeatletter                                                   
\newtheorem*{rep@theorem}{\rep@title}
\newcommand{\newreptheorem}[2]{%
\newenvironment{rep#1}[1]{%
 \def\rep@title{#2 \ref{##1}}%
 \begin{rep@theorem}}%
 {\end{rep@theorem}}}
\makeatother

\newtheorem{thm}{Theorem}[section]
\newtheorem{lemma}[thm]{Lemma}
\crefname{lemma}{Lemma}{Lemmata}
\newtheorem{prop}[thm]{Proposition}

\newtheorem{claim}[thm]{Claim}

\newtheorem*{thm*}{Theorem}
\newtheorem*{lemma*}{Lemma}
\newtheorem*{prop*}{Proposition}
\newtheorem*{corr*}{Corrolary}
\newtheorem*{claim*}{Claim}

\theoremstyle{remark}

\newtheorem*{rmk*}{Remark}
\newtheorem*{conj*}{Conjecture}
\newtheorem*{quest*}{Question}

\theoremstyle{definition}

\newtheorem{dfn}[thm]{Definition}
\newtheorem{exmp}[thm]{Example}
\newtheorem{ex}[thm]{Example}

\newtheorem*{defn*}{Definition}
\newtheorem*{exmp*}{Example}

\newreptheorem{theorem}{Theorem}
\newreptheorem{corollary}{Corollary}
\newreptheorem{proposition}{Proposition}

\newtheorem*{decision*}{Problem}

\newcommand{\R}{\mathbb{R}}

\newcommand{\Z}{\mathbb{Z}}
\newcommand{\N}{\mathbb{N}}

\newcommand{\FPtt}{\mathtt{FP}}

\newcommand{\Rrm}{\mathrm{R}}
\newcommand{\drm}{\mathrm{d}}

\newcommand{\nov}{\widehat{\Z G}^\phi}

\newcommand{\im}{\mathrm{im}}

\newcommand{\supp}{\mathrm{supp}}
\newcommand{\aker}{\mathrm{Aker}}

\newcommand{\Cbf}{\mathbf{C}}
\newcommand{\Fbf}{\mathbf{F}}

\newcommand{\Xcl}{\mathcal{X}}

\newcommand{\col}{\colon}

\def\iff{if and only if }

\title{Quasi-BNS invariants}

\author{Nicolaus Heuer}
\email{nicolaus.heuer.maths@gmail.com}

\author{Dawid Kielak}
\address{
University of Oxford, Oxford, OX2 6GG,
UK}
\email{kielak@maths.ox.ac.uk}

\begin{document}

\begin{abstract}
We introduce the notion of quasi-BNS invariants, where we replace homomorphism to $\R$ by homogenous quasimorphisms to $\R$ in the theory of Bieri--Neumann--Strebel invariants. We prove that the quasi-BNS invariant $Q\Sigma(G)$ of a finitely generated group $G$ is open; we connect it to approximate finite generation of almost kernels of homogenous quasimorphisms; finally we prove a Sikorav-style theorem connecting  $Q\Sigma(G)$ to the vanishing of the suitably defined Novikov homology.
\end{abstract}

\maketitle


\section{Introduction} 

In 1987, Bieri--Neumann--Strebel~\cite{Bierietal1987} introduced the geometric invariant $\Sigma$, nowadays more often called the BNS invariant. The invariant is a set of non-trivial homomorphisms $G \to \R$ with the following fundamental properties (with the last statement proved by Sikorav~\cite{Sikorav1987}). 

\begin{thm}[\cites{Bierietal1987,Sikorav1987}]
\label{BNS}
Let $G$ be a finitely generated group and let $\phi \colon G \to \R$ be a non-zero homomorphism.
\begin{enumerate}
  \item \label{BNS open} The BNS invariant  $\Sigma(G)$ is open.
 \item \label{BNS main point} If $\im \phi = \Z$ then $\phi$ is an algebraic fibring (that is, $\ker \phi$ is finitely generated) \iff $\{\phi, -\phi\} \subseteq \Sigma(G)$.
 \item \label{BNS Sikorav} We have $\phi \in \Sigma(G)$ \iff the first Novikov homology of $(G,\phi)$ vanishes, that is, $\mathrm{H}_1(G;\nov) = 0$.
\end{enumerate}
\end{thm}

Item \eqref{BNS open} is used to construct new algebraically fibred characters out of an old one (when the abelianisation of $G$ has rank at least $2$).
The combination of items \eqref{BNS main point} and \eqref{BNS Sikorav} gives a practical way of finding algebraic fibrings. This approach was used by Friedl--Tillmann~\cite{FriedlTillmann2020} in the setting of one-relator groups, and by the second author in the setting of RFRS groups \cite{Kielak2020a}. It is worth noting that the inspiration behind both papers lies in the realm of $3$-manifolds -- the former article was inspired by the work of Thurston~\cite{Thurston1986} (which also plays a role in \cite{Bierietal1987}), whereas the latter drew from the work of Agol~\cite{Agol2013}.

\smallskip
Since their inception, BNS invariants have been generalised in two main directions.
First, we might want to understand $\ker \phi$ better.
In \cref{BNS}\eqref{BNS main point}, Bieri--Renz~\cite{BieriRenz1988} introduced the \emph{higher BNS invariants} which are capable of computing higher finiteness properties of $\ker \phi$.
Second, we might want to look beyond homomorphisms to $\R$.
In \cite{BieriGeoghegan2016}, Bieri--Geoghegan noticed that a homomorphism (or \emph{character}) $\phi \colon G \to \R$ can be identified with an action of $G$ by translation on a $\mathrm{CAT}(0)$ space $\R$. They then developed a theory in which the role of $\phi$ is taken by an isometric action of $G$ on a proper $\mathrm{CAT}(0)$ space.

This article seeks to generalise BNS invariants  by allowing greater freedom in the choice of $\phi$ as well; the extra flexibility is found however not by looking at more general geometric actions, as was the case for Bieri--Geoghegan, but by leaving geometry entirely. We will replace homomorphisms $\phi \colon G \to \R$ by \emph{quasimorphism}, which can be roughly described as homomorphisms up to a uniform error. The big advantage quasimorphisms have other homomorphisms to $\R$ (and actions on $\mathrm{CAT}(0)$ spaces) is that they are ubiquitous. Take for example an infinite hyperbolic group $G$ with property $(T)$ --- such groups were constructed by \.Zuk~\cite{Zuk2003} (see also \cite{KotowskiKotowski2013} for another proof), and we now have more explicit examples due to Caprace~\cite{Caprace2019}. The abelianisation of $G$ is finite, and more generally its actions on $\mathrm{CAT}(0)$ spaces are severely restricted \cite{NibloReeves1997}. But --- the group $G$ admits infinitely many non-homothetic quasimorphisms (some of them very explicit), as constructed by Calegari--Fujiwara~\cite{CalegariFujiwara2010}.

Allowing quasimorphisms leads to the notion of \emph{quasi BNS invariants}, denoted $Q\Sigma(G)$. This invariant satisfies the following theorem.

\begin{thm}
\label{QBNS}
Let $G$ be finitely generated and let $\phi \colon G \to \R$ be a homogenous quasimorphism with non-zero defect.
\begin{enumerate}
  \item \label{QBNS open} The set $Q\Sigma(G)$ is an open subset of $\R^G$ endowed with the product topology.
 \item \label{QBNS main point} The approximate kernel $\aker \phi$ of $\phi$ is approximately finitely generated \iff $\{\phi, -\phi\} \subseteq Q\Sigma(G)$.
 \item \label{QBNS Sikorav} We have $\phi \in Q\Sigma(G)$ \iff the first Novikov homology of $(G,\phi)$ vanishes, that is, $\mathrm{H}_1(G;\nov) = 0$.
\end{enumerate}
\end{thm}

There are various notions in the above statement that require introduction. The approximate kernel $\aker \phi$ is a natural generalisation of a kernel, and includes elements of $G$ whose image under $\phi$ is very close to $0$. Approximate finite generation is again a natural generalisation of finite generation, and is defined using connectivity of the Rips complex. Novikov rings and Novikov homology are defined in the same way as in the classical setting.

The above theorem shows very clearly that quasi-BNS invariants are the correct analogues of BNS invariants in the setting of approximate group theory, as developed for example by Cordes--Hartnick--Toni\'c~\cite{Cordesetal2020}. This is an emerging topic with great promise.

\subsection*{Acknowledgements}
The authors would like to thank Matt Zaremsky for comments on an earlier version of this article.

This work has received funding from the European Research Council (ERC) under the European Union's Horizon 2020 research and innovation programme, Grant agreement No. 850930.

\section{Quasimorphisms and approximate groups}
We recall the definition of an approximate group and define finiteness properties of approximate groups in Subsection \ref{subsec:approx_groups:def}.
This paper studies approximate groups which arise as \emph{approximate kernels} of quasimorphisms; see Subsection \ref{subsec: approx groups via quasimorphisms}.

\subsection{Approximate groups} \label{subsec:approx_groups:def}

\begin{dfn}
 Let $G$ be a group. A subset $A \subseteq G$ is called an \emph{approximate group} if there is a finite subset $X  \subseteq G$  such that $A \cdot A \subseteq A \cdot X$.
Here and in what follows, $A \cdot B = \{ a b \mid a \in A, b \in B \}$ for $A, B \subseteq G$.
\end{dfn}

Let $G$ be a group with a finite generating set $S$ (which is not necessarily symmetric), and let $\drm$ denote the word metric on $G$ with respect to $S$.
Let $A \subseteq G$ be a subset. The \emph{Rips complex} $\Rrm_S(A,n)$ of $A$ with parameter $n$ is a simplicial complex defined as follows: its $0$th skeleton consists of the elements of $A$, and 
a subset $\{ x_1, \ldots, x_{k+1} \} \subseteq A$ of pairwise distinct $k+1$ elements spans a simplex \iff its diameter is less than $n$.
Here, the diameter of the set $\{ x_1, \ldots, x_{k+1} \}$ is defined as $\max \{ d(x_i, x_j) \mid 1 \leqslant i,j \leqslant k+1 \}$.

\begin{dfn}
 We say that a subset $A$ of $G$ is \emph{approximately connected} \iff there exists a finite generating set $S$ of $G$ and $N \in \N$ such that the Rips complex $\Rrm_S(A,N)$ is connected.
\end{dfn}
It is easy to see that the Rips complex $\Rrm_S(A,n)$ is connected for every $n \geqslant N$, and hence the definition above is equivalent to saying: for every finite generating set $S$ there exists $N \in \N$ such that $\Rrm_S(A,N)$ is connected.

\begin{dfn}[\cite{Cordesetal2020}*{Definition 3.28}]
We say that an approximate group $A \subseteq G$ of a finitely generated group $G$ is \emph{approximately finitely generated} \iff $A$ is approximately connected.
\end{dfn}
In \cite{Cordesetal2020} this notion is referred to as \emph{geometric} approximate finite generation, but we will use the shorter form here.

\subsection{Paths} \label{subsec:notation:paths}

In this subsection we set up some notation for paths in a group $G$ endowed with a finite generating set $S$.

\begin{dfn}
 A \emph{path} $p$ in $G$ is a sequence $(g_0, g_1, \ldots, g_k)$ of elements of $G$ such that $g_i^{-1} g_{i+1} \in S^\pm = S \cup S^{-1}$ for every $i$. For such a path $p$ we will set $o(p) = g_0$ to be the origin and $t(p) = g_k$ to be the terminus of $p$. We will say that $p$ \emph{connects} $g_0$ to $g_k$. We set $p^{-1} = (g_k, \dots, g_0)$.
 
 Given a set-theoretic function $\alpha \col G \to \R$ we set 
 \[\max \alpha(p) = \max \{ \alpha(g_i) \mid i = 1, \ldots, k \}\]
  and similarly $\min \alpha(p) = \min \{ \alpha(g_i) \mid i = 1, \ldots, k \}$.
 
 Given another path $q = (h_0, \dots, h_n)$ in $G$, we define the \emph{concatenation} of $p$ and $q$ to be
 \[
  pq = (g_0, \dots, g_k, g_k {h_0}^{-1} h_1, \dots, g_k {h_0}^{-1} h_n)
 \]
\end{dfn}

Note that the last definitions applies in particular to the situation in which $p= ( g_0)$ is a single-vertex path, in which case $pq = (g_0)q$ is a translate of $q$ starting at $g_0$. Also, the last definition allows us to write $p^n$ for every $n \in \Z \smallsetminus \{0\}$.

We can reformulate the notion of approximate connectivity using paths (we omit the obvious proof).
\begin{lemma}
 Let $A$ be  a subset of a group $G$ endowed with a finite generating set $S$. The subset $A$ is approximately connected \iff there exists a constant $K \geqslant 0$ such that for every $a, b \in A$ there exists a path $p$ connecting $a$ to be $b$ such that for every vertex $v$ on the path we have $d(v,A) \leqslant K$.
\end{lemma}

\subsection{Quasimorphisms and approximate kernels} \label{subsec: approx groups via quasimorphisms}
\begin{dfn}
A map $\phi \col G \to \R$ is called a \emph{quasimorphism} if there is some $C > 0$ such that
for all $g,h \in G$, $|\phi(g) + \phi(h) - \phi(gh)| \leqslant C$.
We call the smallest such $C$ the \emph{defect of $\phi
$} and denote it by $D(\phi)$.
We say that $\phi$ is homogeneous if $\phi(g^n) = n \phi(g)$ for every $g \in G$, $n \in \Z$.
\end{dfn}

Every quasimorphism $\phi$ is in bounded pointwise distance to a unique homogeneous quasimorphism: the function $\bar{\phi}$ defined pointwise via
\[
\bar{\phi}(g) = \lim_{n \to \infty} \frac{\phi(g^n)}{n}
\]
is a quasimorphism of defect at most $2 D(\phi)$ and $| \phi(g) - \bar{\phi}(g) | \leqslant D(\phi)$ for every $g \in G$, cf. \cite{scl}*{Lemma 2.21} of Calegari.
The set of homogeneous quasimorphisms forms a vector space under pointwise addition and will be denoted by $Q^h(G)$.
We will topologise $Q^h(G)$ in a non-standard way. We write it as a disjoint union
\[
 Q^h(G) =\bigsqcup_{x \in [0,\infty)} \{ \phi \in Q^h(G) \mid D(\phi) = x \} 
\]
and topologise each subset $\{ \phi \in Q^h(G) \mid D(\phi) = x \} $ by embedding it into $\R^G$ with product topology and then using subspace topology. Note that this in particular separates homomorphisms from proper quasimorphisms.

\begin{lemma}[\cite{scl}*{Lemma 2.24}] \label{lemma:phi large on commutators}
For a homogeneous quasimorphism $\phi$ on a group $G$ we have
\[
D(\phi) = \sup_{g,h \in G} \phi([g,h])
\]
\end{lemma}

\begin{prop} \label{prop:approx kernel is approx subgrp} 
For every $\phi \in Q^h(G)$ the set
\[
\aker(\phi) = \big\{ g \in G \, \big|\, |\phi(g)| \leqslant 2 D(\phi) \big\}
\]
is an approximate subgroup. We will call this set the \emph{approximate kernel} of $\phi$.
\end{prop}

\begin{proof}
Suppose that $D(\phi) = 0$. Then $\phi$ is a homomorphism and hence $\aker(\phi)$ is the honest kernel of $\phi$. 

Now suppose that $D(\phi) \neq 0$.
By \cref{lemma:phi large on commutators}, there exists $c = [g,h]$ such that 
\[\frac{ 4 D(\phi)}{5} < \phi(c) \leqslant D(\phi).\]
We claim that $\aker(\phi) \cdot \aker(\phi) \subseteq \aker(\phi) \cdot X$ where 
\[X = \{c^5,  c^4, \ldots,  c^{-4}, c^{-5} \}\]
Let $g, h \in \aker(\phi)$. The definitions tell us that $\phi(g h c^n) = \phi(g h) + n \phi(c) + d_1{(n)}$ and $\phi(g h) = \phi(g) + \phi(h) + d_2$, where $d_1{(n)}$ depends on $n$, and where $|d_1{(n)}|, |d_2| \leqslant D(\phi)$ for every $n$. Hence
\begin{align*}
|\phi(g h c^n) | &= |\phi(g) + \phi(h) + n \phi(c)+ d_1{(n)} + d_2| \\ &\leqslant |\phi(g) + \phi(h) + n \phi(c)+ d_2| + D(\phi)
\end{align*}
Also, as $|\phi(g) + \phi(h) + d_1|\leqslant 5D(\phi)$, we may choose $n \in \{ -5, \ldots, 5 \}$ such that 
\[
 |\phi(g) + \phi(h) + n \phi(c)+ d_2| \leqslant D(\phi)
\]
and so 
$|\phi(g h c^n)| \leqslant 2 D(\phi)$. Hence $g h c^n \in A K (\phi)$ and thus $gh \in \aker(\phi) \cdot X$.
\end{proof} 

Interestingly, for our purposes the size of $X$ in the definition of approximate groups could have been bounded uniformly by $11$.

Observe that the approximate kernel of $\phi$ is the same as the approximate kernel of $\lambda \phi$ for every $\lambda \in \R \setminus \{ 0 \}$.

For later purposes the following characterisation of approximate finite generation will be useful:
\begin{prop} \label{prop:finiteness properties of approx kernels}
Let $\phi \col G \to \R$ be a homogeneous quasimorphism with non-zero defect $D = D(\phi)$, and suppose that $G$ is equipped with a finite generating set $S$. Then $\aker(\phi)$ is approximately finitely generated if and only if there is a constant $K > 0$ such that for every two elements $g, h \in \aker(\phi)$ there is a path $p$ in $G$ from $g$ to $h$ such that $\min \phi(p) \geqslant -K$ and $\max \phi(p) \leqslant K$.
\end{prop}
\begin{proof}
As in the proof of \cref{prop:approx kernel is approx subgrp}, there exists an element $c \in G$ such that
\[\frac{ 4 D}{5} < \phi(c) \leqslant D.\]
Let $C = \drm (1,c)$.

Suppose first that we have a constant $K$ as in the statement. Let $g,h \in \aker(\phi)$ be given. The assumption gives us a path $p$ connecting $g$ to $h$ in $G$ such that $\min \phi(p) \geqslant -K$ and $\max \phi(p) \leqslant K$. Let $k$ be any vertex of $p$. There exists $m \in \Z$ with $|m| \leqslant \frac {5K} {4D}$ such that
\[
 |\phi(kc^m)| \leqslant 2D
\]
and therefore $kc^m \in \aker(\phi)$. Note that $\drm(k,kc^m) \leqslant |m|C \leqslant \frac{5KC} {4D}$. This means that for every $k \in p$ we may find $k' \in \aker(\phi)$ such that $\drm(k,k') \leqslant \frac{5KC} {4D}$. If $k_1$ and $k_2$ are successive points on $p$, then $\drm(k_1',k_2') \leqslant \frac{5KC} {2D} + 1$, and so 
$g$ and $h$ lie in the same connected component of $\Rrm_S(\aker(\phi),n)$ for every $n \geqslant \frac{5KC} {2D} + 1$, as required.

\smallskip
Now suppose that $\Rrm_S(\aker(\phi),n)$ is connected for some $n$. Let 
\[|\phi(S)| = \max_{s \in S} |\phi(s)|\]
 and $K = n(|\phi(S)| + D) + 2D$. Take $g , h \in \aker(\phi)$. Since $\Rrm_S(\aker(\phi),n)$ is connected, there exists a path $q$ in $\Rrm_S(\aker(\phi),n)$ connecting $g$ to $h$. This path can easily be transformed to a path in $G$ connecting $g$ to $h$, where each edge in $q$ is replaced by a path of length at most $n$ in $G$. It is immediate that for every vertex $k$ of $p$ we have $|\phi(k)| \leqslant n(|\phi(S)| + D) + 2D = K$.
\end{proof}

The above statement is also true for a homomorphism $\phi \colon G \to \Z$. Indeed, in the proof we take $c \in G$ to be such that $\phi(c)$ generates $\Z$, and otherwise argue in a similar vain. Interestingly though the statement is false for homomorphisms with bigger image, as we will now show.

\begin{ex}
 Let $G = F_2 \times \Z$, where $F_2 = \langle a,b \rangle$ and $\Z = \langle c \rangle$. Let $\phi \colon G \to \R$ be a homomorphism induced by $\phi(a) =1, \phi(b) = 0, \phi(c) = \sqrt{2}$. Clearly $\ker \phi = \langle a^{-i}ba^i : i \in \Z \rangle$ is not finitely generated. But given any two elements $g,h \in \ker \phi$ and a path $p$ in $F_2 \times \{0\}$ with generators $\{a^{\pm 1}, b^{\pm 1}\}$ connecting $g$ to $h$, we may modify $p$ using the generator $c$ so that the resulting path $q$ satisfies $\min \phi(q) \geqslant -3$ and $\max \phi(q) \leqslant 3$. 
\end{ex}

\section{Quasi-BNS invariants}

\subsection{Definition}

For a homogeneous quasimorphism $\phi \in Q^h(G)$ we let 
\[G_\phi = \{ g \in G \mid \phi(g) > 0 \}\] 
Analogously to the classical definition we set
\[
Q\Sigma(G) = \{ \phi \in Q(G) \mid G_\phi \text{ is approximately connected} \} \subseteq Q^h(G)
\]
Note that this set is invariant under positive scalar multiplication, i.e.\ for every $\lambda \in \R^+$, $\phi \in Q\Sigma(G)$ if and only $\lambda \phi \in Q\Sigma(G)$. Note also that this definition extends the classical one, that is, a non-trivial homomorphism $\phi \colon G \to \R$ lies in $Q\Sigma(G)$ \iff it lies in $\Sigma(G)$.

The following is immediate but useful:

\begin{prop} \label{prop:positive path}
Let $G=\langle S \rangle$ with $S$ finite,  and let $\phi \col G \to \R$ be a homogenous quasimorphism. Then $\phi \in Q \Sigma(G)$ if and only if there is a constant $K > 0$ such that for every two points $g, h \in G_\phi$ there is a path $p$ from $g$ to $h$ such that $\min \phi(p) \geqslant -K$.
In this case, there is a constant $K'>0$ such that for any two elements $g,h \in G$ there is a path $p$ from $g$ to $h$ such that $\min \phi(p) \geqslant \min \{ \phi(g), \phi(h) \} - K'$.
\end{prop}

\subsection{Approximate finite generation of approximate kernels}
We will relate approximate finite generation of the approximate kernel $\aker(\phi)$ to $Q \Sigma(G)$, in analogy with the situation for the classical BNS invariant.
\begin{thm} 
Let $\phi$ be a homogeneous quasimorphism with positive defect. Then $\phi, -\phi \in Q \Sigma(G)$ if and only if the approximate kernel $\aker(\phi)$ is approximately finitely generated.
\end{thm}
\begin{proof}
Suppose that $\phi$ is a homogeneous quasimorphism with $D = D(\phi) > 0$. Let $c = [g,h]$ be a commutator such that $\frac 4 5 D < \phi(c) \leqslant D$, as before.

Assume that $\aker(\phi)$ is approximately finitely generated and let $g_1, g_2 \in G_\phi$. Analogously to the proof of  \cref{prop:approx kernel is approx subgrp} there are $n_1, n_2 \geqslant 0$ such that 
\[g_1', g_2' \in \aker(\phi)\]
 with $g_i' = g_i c^{-n_i}$.
Let $p_c$ be a path from $1$ to $c$.
By \cref{prop:finiteness properties of approx kernels} there is a uniform $K > 0$ such that we can find a path $p$ from $g_1'$ to $g_2'$ with $\min \phi(p) \geqslant -K$.
Finally, the path $p' = (g_1)\cdot p_c^{-n_1} \cdot p \cdot p_c^{n_2}$ connects $g_1$ to $g_2$, and we have $\min \phi(p')$  uniformly bounded. By \cref{prop:positive path}, we have shown that $\phi \in Q \Sigma(G)$; by symmetry $-\phi \in Q \Sigma(G)$ as well.

\smallskip
Now assume that both $\phi, -\phi \in Q\Sigma(G)$.
To prove that $\aker \phi$ is approximately finitely generated, we will use a peak reduction argument.

We take a finite symmetric generating set $S$ for $G$ containing $c$. Given a path,
we define its vertex $v$ to be \emph{inessential} \iff
it is preceded by $vc^{\pm 1}$ and followed by $vc^{\pm 1}$ where the signs can be different. A vertex that is not inessential will be called \emph{essential}.

Since $\phi, -\phi \in Q\Sigma(G)$, let $K'$ be the second constant from \cref{prop:positive path}
which works for both $\phi$ and $-\phi$; for convenience we take $K' > 2D$. Let 
\[
n = \Big\lfloor \frac { 5} {4 D}\big(K' +\max_{s,t \in S} \phi(st) + D\big)\Big\rfloor + 3
 \]

Since $-\phi \in Q\Sigma(G)$, for every pair $(s,t) \in S^2$ there exists a path $q'_{s,t}$ such that 
\[
q_{s,t} = (1, c^{-1}, \dots, c^{-n}) q'_{s,t} (1, c, \dots, c^n)
\]
connects $1$ to $st$, and its essential vertices $v$ satisfy 
\begin{align*}
\phi(v) &\leqslant \max\{\phi(c^{-n}), \phi(stc^{-n})\} +K' \\
&\leqslant \max_{s,t \in S} \phi(st) - n\phi(c) +D +K' \\
&\leqslant \max_{s,t \in S} \phi(st) - \frac {4nD} 5 +D+K' \\
&\leqslant \max_{s,t \in S} \phi(st) - \frac {4D} 5\Big(\frac { 5} {4 D}\big(K' +\max_{s,t \in S} \phi(st) + D\big) + 2\Big) +K' \\
&< -\frac 8 5 D \\ 
&< -D
 \end{align*}
Let $N$ be such that $-N < \min \phi( q_{s,t})$ for every $(s,t)$ and such that $N$ is greater than $K' + 2D$.

Now take $x,y \in \aker(\phi)$. We claim that there exists a path $p$ connecting $x$ to $y$ such that
\[
 -N  < \min \phi(p) \leqslant  \max \phi(p)  \leqslant M + 2D
\]
where $M = 3 D + \max_{s \in S} |\phi(s)|$.

To prove the claim, we need two definitions: Given a path $p'$ connecting $x$ to $y$, the
\emph{height} of $p'$ is defined to be $\max \lfloor \phi(v) \rfloor$, with the maximum ranging over essential vertices. An essential vertex $v$ of $p'$ is a \emph{peak} \iff $\lfloor \phi(v) \rfloor$ coincides with the height of $p'$. 

Let $p$ be a path from $x$ to $y$ with
\[
-N  < \min \phi(p)                                                                                                                                                                                                                                                                                                                                                                                                                                                                                                                                                                                                                                                                                                                                                                                                                         \]
and such that $p$ has the smallest possible height, and then the smallest possible number of peaks among all paths satisfying the inequality above; the path $p$ exists since $-N <  -2 D -K'  \leqslant \min \{\phi(x), \phi(y) \}  -K'$.
Suppose that the height of $p$ is greater than $M$. Let $v_1$ be the first peak of $p$, and note that $v_1$ is neither $x$ nor $y$, since $|\phi(x)|$ and $|\phi(y)|$ are strictly smaller than $3 D$. Let $v_0$ be the vertex preceding $v_1$, and let $v_2$ be the vertex following $v_1$ so we have $v_2 = v_1t$ and $v_1 = v_0 s$ for some generators $s$ and $t$. Now we replace the subpath $(v_0, v_1, v_2)$ in $p$ by 
$(v_0) q_{s,t}$, and form a new path $p'$ this way. By construction, the essential vertices of $(v_0) q_{s,t}$ cannot be peaks, since we have already shown that an essential vertex $w$ of $q_{s,t}$ satisfies $\phi(w) < -D$, and so the essential vertices of $(v_0) q_{s,t}$ have to have $\phi$-values strictly lower than $\phi(v_0)$, and $v_0$ is not a peak.
Also, every essential vertex of $p'$ is either an essential vertex of $p$, or an essential vertex of the subpath $(v_0) q_{s,t}$. Hence $p'$ has fewer peaks or lower height than $p$. Also, 
\[
\min \phi(p')  = \min \Big\{ \min \phi(p), \min \phi\big( (v_0) q_{s,t} \big) \Big\}                                                                                                                                                                                                                                                                                                                                                                                                                                                                                                                                                                                                                                                                                                                                       \]
and we have
\[
\min \phi(p) > -N
\]
Let us now compute the second minimum:
\begin{align*}
\min \phi\big( (v_0) q_{s,t} \big) &\geqslant \phi(v_0) + \min \phi(q_{s,t}) - D\\
&> \phi(v_1) - \max_{s \in S} |\phi(s)| - D  - N - D \\
&> M -  \max_{s \in S} |\phi(s)|   - N -2D \\
&> D-N \\
&>-N                             
\end{align*}
Together these inequalities yield $\min \phi(p') > -N$, contradicting the minimality of $p$. We conclude that the height of $p$ is at most $M$. 

We now remove all inessential vertices from $p$ that appear when the path $p$ backtracks along an edge corresponding to $c$. Note that a vertex of the new path is essential if and only if at least one of the vertices it came from was essential.
Let $v$ be a vertex of $p$ such that $p(v) > M$. Since $v$ is inessential, it lies on a segment using only the generators $c$ or $c^{-1}$ connecting two essential vertices, say $u$ and $w$. Without loss of generality, let us take $v = uc^\alpha$ and $w = uc^\beta$, with $0 \leqslant \alpha \leqslant \beta$. We have
\[
M \geqslant \phi(w) \geqslant \phi(u) + \beta \phi(c) - D 
\]
and so 
\[
M - \phi(u) +D \geqslant  \beta \phi(c) \geqslant \alpha \phi(c) \geqslant \phi(v) - \phi(u) - D 
\]
yielding  $\phi(v) \leqslant M + 2D$. This finishes the proof.
\end{proof}
As a corollary we obtain \cref{QBNS}\eqref{QBNS main point}.

\section{Sikorav's Theorem for quasi BNS invariants}
The aim of this section is to generalise  the following theorem due to Sikorav~\cite{Sikorav1987} which provides a homological characterisation of elements $\phi \in \Sigma(G)$:
\begin{thm}{(Sikorav~\cite{Sikorav1987})}
Let $\phi \col G \to \R$ be a non-trivial homomorphism, where $G$ is a finitely generated group. Then $\phi \in \Sigma(G)$ if and only if $H_1(G, \nov) = 0$.
\end{thm}

Here $\nov$ denotes the Novikov ring with respect to $\phi$ defined as follows: the elements of $\nov$ are formal sums $\sum_{g \in G} \lambda_g g$ with $\lambda_g \in \Z$ and such that for every $x \in \R$ the set
\[
 \{ g \in G \mid \phi(g) < x, \lambda_g \neq 0\}
\]
is finite. Addition is pointwise, and multiplication extends the group ring multiplication of $\Z G$. Note that the exact same definition gives a ring when $\phi$ is a quasimorphism.

We will generalise Sikorav's result to the case of quasi-BNS invariants.
\begin{thm} \label{thm:quasi-sikorav}
Let $\phi \col G \to \R$ be a non-trivial homogenous quasimorphism, where $G$ is a finitely generated group. The following are equivalent:
\begin{enumerate}
 \item $\phi \in Q \Sigma(G)$
 \item There exists an open neighbourhood $U$ of $\phi$ such that for every $\psi \in U$ we have $\mathrm{H}_1(G, \widehat{\Z G}^\psi) = 0$.
 \item $\mathrm{H}_1(G, \nov) = 0$.
\end{enumerate}
\end{thm}
\begin{proof}
We take a finite symmetric generating set $S$ of $G$ which contains an element $c$ with $\phi(c) > 0$.

Let $X^{(1)}$ denote the Cayley graph of $G$ with respect to $S$, modified slightly so that if $s \in S$ is an involution, then the edges corresponding to $s$ and $s^{-1}$ emanating from the same vertex are identified. We extend $X^{(1)}$ to a 2-dimensional acyclic $G$-CW-complex $X$ by adding $G$-orbits of $2$-cells.
We identify the $0$-skeleton $X^{(0)}$ of $X$ with $G$, and this way we can view $\phi$ as a map $X^{(0)} \to \R$. We extend $\phi$ to $X$ inductively by defining $\phi$ of a cell $C$ to be the minimum of the value of $\phi$ on all the $0$-cells contained in $C$.

This choice of $X$ allows us also to continue using our notation for paths, since any two vertices in $X^{(0)} = G$ differing by a generator form $S$ are connected in $X^{(1)}$ by a unique edge.

We now define a \emph{Novikov $n$-chain} in $X$ to be a formal linear combination of $n$-cells in $X$ with coefficients in $\Z$ such that for every $x \in \R$ there are only finitely many $n$-cells in the chain with non-zero coefficient that are mapped to $(-\infty,x)$ by $\phi$. It is immediate that $\mathrm{H}_1(G, \nov)$ coincides with the space of Novikov $1$-cycles divided by Novikov $2$-boundaries, with cycles and boundaries defined in the usual way.
We are now ready to start the proof.

\smallskip
\noindent {$\mathbf{(1)\Rightarrow (2)}$}
Suppose that $\phi \in Q \Sigma(G)$. Let $K\geqslant 0$ be such that for every two elements of $a,b \in G_\phi$ there exists a path $p$ connecting them such that $\min \phi(p) \geqslant \min\{ \phi(a), \phi(b) \} -K$. Pick $n$ such that $n\phi(c) > K+D(\phi) + 1$. Let $U$ be an open neighbourhood of $\phi$ defined by
\[
 U = \{ \psi \in Q^h(G) \mid D(\psi) = D(\phi), n\psi(c) > K+D(\psi) + 1 \}
\]

For every  $s \in S$ we construct a $1$-cycle $z_s$ by combining two paths connecting $c^n$ to $sc^n$. The first path is $(c^n)(1,c^{-1})^{n}(1,s)(1,c)^n$; the second path $p_s$ comes from the fact that $\phi \in Q\Sigma(G)$, and satisfies $\min \phi(p_s) \geqslant n\phi(c) - K >  D(\phi) + 1$. We intersect $U$ with  further open sets and obtain $V \subseteq U$ such that for all $\psi \in V$ and $s \in S$ we have
\[
 \min \psi(p_s) >  D(\psi) + 1
\]
Since $z_s$ is a $1$-cycle in $X$, and $X$ has vanishing first homology, there exists a $2$-chain $y_s$ in $X$ with $\partial y_s = z_s$.

Now every edge in $p_s$ is of the form $g(1,s')$, $s' \in S$ and $g \in G$ with $\psi(g) >  D(\psi) +1$ for every $\psi \in V$. For every such edge we add $-gy_{s'}$ to $y_s$ and form a new $2$-chain whose boundary is a $1$-cycle made of the path $(c^{2n})(1,c^{-1})^{2n}(1,s)(1,c)^{2n}$ and a path whose vertices have $\psi$-values bounded below by $2\big(D(\psi) + 1\big) - D(\psi) = D(\psi) + 2$. Moreover, $\min \psi(gy_{s'}) \geqslant \psi(g) + \min \psi(y_{s'}) - D(\psi) > \min \psi(y_{s'}) + 1$. Proceeding inductively, we build a Novikov $2$-chain $y'_s$ whose boundary is the Novikov $1$-cycle $x'_s$ obtained from the edge $(1,s)$ and infinite rays emanating from $1$ and $s$ and using only the edges corresponding to $c$. The chains are Novikov with respect to every $\psi \in V$.

Now fix $\psi \in V$ and let $z$ be any $\psi$-Novikov $1$-cycle in $X$. Adding the boundaries of Novikov $2$-chains $y'_s$ with $s \neq c$ we just constructed we easily arrange for $z$ to consist solely of edges corresponding to $c$. But the only Novikov $1$-cycle of this form is the $0$ cycle.

\smallskip
\noindent {$\mathbf{(2)\Rightarrow (3)}$} This is obvious.

\smallskip
\noindent {$\mathbf{(3)\Rightarrow (1)}$}
Now suppose that $\mathrm{H}_1(G, \nov) = 0$, that is, that every Novikov $1$-cycle is the boundary of a Novikov $2$-chain. Let $a, b \in G_\phi$ be two elements, and let $q$ be any path connecting them in $X$. Consider the Novikov $1$-cycle $z$ formed by $q$ and the two infinite rays emanating from $a$ and $b$ and using only the edges corresponding to $c$. By assumption, there exists a Novikov $2$-chain $y'$ such that $\partial y' = z$.

Let $y$ be a $2$-chain obtained from $y'$ by keeping only the cells whose value under $\phi$ is strictly negative. Note  that $\partial y - z$ is supported solely on cells with $\phi$-value positive or zero. Now the support of $\partial y - z$ must contain a path $p'$ connecting $ac^m$ to $bc^n$ for some non-negative $m$ and $n$. We construct the path $p$ connecting $a$ to $b$ by first following the ray from $a$ to $ac^m$, then by following $p'$, and finally by following the ray in the negative direction from $bc^n$ to $b$. It is clear from the construction that $\min \phi(p) \geqslant -D(\phi)$. This finishes the proof.
\end{proof}

As a corollary we obtain \cref{QBNS}\eqref{QBNS open} and \eqref{QBNS Sikorav}.

\section{Examples}

\begin{ex}
 Let $F_n$ be the free group of rank $n \geqslant 2$. The quasi-BNS invariant $Q\Sigma(F_n)$ is empty. This can be seen as follows: Let $\phi \in Q^h(F_n)$ be a non-trivial homogeneous quasimorphism, and  take $c \in F_n$  such that $\phi(c) > 0$. Using \cref{lemma:phi large on commutators} we see that $\aker (\phi)$ contains the commutator subgroup of $F_n$. For every element $x \in \aker(\phi)$ not commuting with $c$ we have that $c^{-n} x c^n$ lies in the $2D(\phi)$-neighbourhood of $\aker(\phi)$. But there is no path from $x$ to $c^{-n} x c^n$ which stays uniformly close to $\aker(\phi)$ for all $n \in \N$, since the natural Cayley graph of $F_2$ is a tree.
\end{ex}

\begin{ex}
 For amenable groups every homogeneous quasimorphism is a homomorphism, hence in this case $Q \Sigma(G) = \Sigma(G)$.
\end{ex}

\bibliography{QBNS}

\end{document}